\documentclass[11pt]{article}
\usepackage[francais]{babel}
\usepackage[latin1]{inputenc}
\usepackage[T1]{fontenc}
\usepackage{a4,amsmath,amsfonts}
\title{Un lemme matriciel effectif}
\author{Pascal Autissier}
\begin{document}

\maketitle

\newcommand{\D}{\displaystyle}

{\bf Abstract:} We show an almost optimal effective version of Masser's
matrix lemma \cite{Mass}, giving a lower bound of the height of an abelian
variety in terms of its period lattices.\\

{\bf Résumé:} On donne ici une version effective presque optimale du lemme
matriciel de Masser, qui consiste à minorer la hauteur d'une variété abélienne
en fonction de ses réseaux des périodes.\\


\section{Introduction}

Soit $K$ un corps de nombres. Désignons par $G_K$ l'ensemble des plongements
$\sigma$ de $K$ dans $\mathbb{C}$. Masser a montré dans \cite{Mass} une
minoration, connue sous le nom de lemme matriciel, de la hauteur d'une variété
abélienne $A$ sur $K$ en termes de matrices des périodes des $A_\sigma$. Ce
résultat est l'un des ingrédients utilisés par Masser et Wüstholz \cite{MaWu}
pour prouver leur fameux théorème des périodes.\\

Bost \cite{Bo1} a revisité ce lemme matriciel en minorant la hauteur de
Faltings (stable), notée $h_{\rm Fa}(A)$, de $A$ en fonction des diamètres
d'injectivité des $A_\sigma$. Graftieaux \cite{Graf}, David-Philippon
\cite{DaPh} et Gaudron \cite{Gaud} en ont ensuite donné des versions effectives
({\it cf.} remarque 1.3 ci-dessous).\\

On se propose ici de raffiner ces versions de manière asymptotiquement
optimale. Introduisons d'abord quelques notations:

Dans tout la suite, on désigne par $\kappa$ la constante
$\D\kappa=\sqrt{\frac{3}{2\pi^3e}}$. Lorsque $(A;L)$ est une variété abélienne
complexe principalement polarisée, on note $\rho(A;L)$ son diamètre
d'injectivité ({\it i.e.} le premier minimum de son réseau des périodes).\\

{\bf Théorème 1.1:} {\it Soit $(A;L)$ une $K$-variété abélienne de dimension
$g\ge1$, principalement polarisée. En posant
$\D\rho_\sigma=\min\Bigl(\rho(A_\sigma;L_\sigma);\sqrt{\frac{\pi}{3g}}\Bigr)$ pour
tout $\sigma\in G_K$, on a alors l'inégalité suivante:
$$h_{\rm Fa}(A)\ge\frac{1}{[K:\mathbb{Q}]}\sum_{\sigma\in G_K}\Bigl(\frac{\pi}{6\rho_\sigma^2}+g\ln(\kappa\rho_\sigma\sqrt{g})\Bigr)\quad.$$}\\

Ce résultat est utilisé par Gaudron et Rémond \cite{GaRe} pour donner de
nouvelles versions effectives du théorème des périodes.\\

{\bf Remarque 1.2:} La constante $\D\frac{\pi}{6}$ devant les
$\D\frac{1}{\rho_\sigma^2}$ est optimale. En effet, fixons un corps de nombres
$K_0$ et une $K_0$-variété abélienne $(A_0;L_0)$ de dimension $g-1$,
principalement polarisée. Alors pour toute extension finie $K$ de $K_0$ et
toute $K$-courbe elliptique $E$ à potentiellement bonne réduction partout, on
a, en posant $A=A_{0K}\times E$, l'estimation
$$h_{\rm Fa}(A)=\frac{1}{[K:\mathbb{Q}]}\sum_{\sigma\in G_K}\Bigl(\frac{\pi}{6\rho_\sigma^2}+\ln\rho_\sigma\Bigr)+O(1)\quad,$$
où le $O(1)$ ne dépend que de $(g;K_0;A_0;L_0)$ (mais pas de $(K;E)$). On le
voit en appliquant le théorème 7.b de \cite{Falt}.\\

{\bf Remarque 1.3:} \`A titre de comparaison, Graftieaux et Gaudron ont obtenu
l'énoncé 1.1 avec, au lieu de ce $\D\frac{\pi}{6}$, une fonction $c(g)$ de $g$
qui converge vers 0 (plus vite que $\D\frac{1}{g^g}$) lorsque $g$ tend vers
$+\infty$; David-Philippon ont trouvé une minoration de la forme
$\D\frac{c_0}{g^4[K:\mathbb{Q}]}\max_{\sigma\in G_K}\frac{1}{\rho_\sigma^2}$ (où la
constante $c_0$ est absolue).\\

Enfin, citons pour mémoire la forme simplifiée (affaiblie mais plus maniable)
suivante:\\

{\bf Corollaire 1.4:} {\it Soit $\varepsilon\in]0;1[$. Soit $(A;L)$ une
$K$-variété abélienne de dimension $g\ge1$, principalement polarisée. On a
alors
$$h_{\rm Fa}(A)+\frac{g}{2}\ln\Bigl(\frac{2\pi^2}{\varepsilon}\Bigr)\ge\frac{(1-\varepsilon)\pi}{6[K:\mathbb{Q}]}\sum_{\sigma\in G_K}\frac{1}{\rho(A_\sigma;L_\sigma)^2}\quad.$$}\\

{\it Démonstration:} \`A partir du théorème 1.1, il suffit d'écrire
$\ln u_\sigma\le u_\sigma-1$ avec
$\D u_\sigma=\frac{\varepsilon\pi}{3g\rho_\sigma^2}$ pour tout
$\sigma\in G_K$. $\square$\\

Après des préliminaires sur les réseaux (section 2), on prouve le théorème 1.1
à la section 3.2.\\

Je remercie \'Eric Gaudron et Gaël Rémond de m'avoir incité à rédiger ce texte.
Je remercie également Fabien Pazuki pour d'intéressantes discussions concernant
le lemme matriciel.\\

\section{Géométrie des nombres}

Soit $g$ un entier $\ge1$. On note $\mathbb{S}_g$ l'ensemble des matrices
$Y\in{\rm M}_g(\mathbb{R})$ symétriques définies positives.\\

{\bf Définitions:} Soit $Y\in\mathbb{S}_g$. Munissons $\mathbb{R}^g$ de la
norme euclidienne $|\!|\ |\!|_Y$ telle que $|\!|x|\!|_Y^2={}^{\rm t}\!xYx$ pour
tout $x\in\mathbb{R}^g$, et posons
$\D\psi_Y(x)=\min_{m\in\mathbb{Z}^g}|\!|x-m|\!|_Y$. Le {\bf premier minimum}
$\lambda_1(Y)$ et le {\bf minimum inhomogène} $\mu(Y)$ de $Y$ sont les réels
définis par
$$\lambda_1(Y)=\min_{m\in\mathbb{Z}^g-\{0\}}|\!|m|\!|_Y\quad\mbox{et}\quad\mu(Y)=\max_{x\in\mathbb{R}^g}\psi_Y(x)\quad.$$

{\bf Lemme 2.1:} {\it Soit $Y\in\mathbb{S}_g$. On a l'inégalité
$2\mu(Y)\lambda_1(Y^{-1})\ge1$.}\\

{\it Démonstration:} Choisissons un $\gamma\in\mathbb{Z}^g$ tel que
$|\!|\gamma|\!|_{Y^{-1}}=\lambda_1(Y^{-1})$. Les coordonnées de $\gamma$ sont
premières entre elles (par minimalité), donc il existe $m\in\mathbb{Z}^g$
vérifiant une relation de Bézout $ ^{\rm t}\!\gamma m=1$. Posons
$\D x=\frac{m}{2}$ et montrons que $2\psi_Y(x)\lambda_1(Y^{-1})\ge1$. Pour tout
$n\in\mathbb{Z}^g$, on a
$$1\le|1-2 ^{\rm t}\!\gamma n|=2| ^{\rm t}\!\gamma(x-n)|\le2|\!|\gamma|\!|_{Y^{-1}}|\!|x-n|\!|_Y=2\lambda_1(Y^{-1})|\!|x-n|\!|_Y\quad.$$

D'où le résultat. $\square$\\

Dans la suite de cette section, on fixe un $Y\in\mathbb{S}_g$ et on note
simplement $|\!|\ |\!|$ et $\psi$ au lieu de $|\!|\ |\!|_Y$ et $\psi_Y$.
Désignons par $\nu$ la mesure de Lebesgue sur $\mathbb{R}^g$, et posons
$F=[0;1]^g$.\\

{\bf Lemme 2.2:} {\it Soit $Y\in\mathbb{S}_g$. On a la minoration suivante:
$$\int_F\psi(x)^2{\rm d}\nu(x)\ge\frac{\mu(Y)^2}{3}\quad.$$}\\

{\it Démonstration:} On prend un $y\in\mathbb{R}^g$ vérifiant $\mu(Y)=\psi(y)$.
Soit $x\in\mathbb{R}^g$. Pour $m$ et $n$ dans $\mathbb{Z}^g$, l'identité du
parallélogramme donne
$$\begin{array}{rcl}
2|\!|x-m|\!|^2+2|\!|x-y-n|\!|^2&=&|\!|y-m+n|\!|^2+|\!|2x-y-m-n|\!|^2\\
&\ge&\psi(y)^2+\psi(2x-y)^2\quad.\\
\end{array}$$

Il en découle $2\psi(x)^2+2\psi(x-y)^2\ge\mu(Y)^2+\psi(2x-y)^2$. On conclut en
intégrant cette inégalité contre $\nu$ et en utilisant la
$\mathbb{Z}^g$-périodicité de $\psi$. $\square$\\

Pour tout $(t;x)\in\mathbb{R}^*_+\times\mathbb{R}^g$, on pose
$$f_Y(t;x)=\sqrt{\det(Y)}\sum_{m\in\mathbb{Z}^g}\exp(-\pi t|\!|x-m|\!|^2)\quad.$$

{\bf Lemme 2.3:} {\it On a les propriétés suivantes.

$(\alpha)$ Soit $x\in\mathbb{R}^g$. L'application
$\mathbb{R}^*_+\rightarrow\mathbb{R}$ qui à $t$ associe
$f_Y(t;x)\exp(\pi t\psi(x)^2)$ est décroissante.

$(\beta)$ Soit $t\in\mathbb{R}^*_+$. On a l'estimation
$\D\int_F\ln f_Y(t;x){\rm d}\nu(x)\le-\frac{g}{2}\ln t$.}\\

{\it Démonstration:} $(\alpha)$ Cette fonction est une somme d'exponentielles
décroissantes.\\

$(\beta)$ En utilisant l'inégalité de convexité de Jensen, on trouve
$$\begin{array}{rcl}
\D\int_F\ln f_Y(t;x){\rm d}\nu(x)&\le&\D\ln\int_Ff_Y(t;x){\rm d}\nu(x)\\
&=&\D\ln\Bigl[\sqrt{\det(Y)}\int_{\mathbb{R}^g}\exp(-\pi t|\!|x|\!|^2){\rm d}\nu(x)\Bigr]=-\frac{g}{2}\ln t\quad.\ \square\\
\end{array}$$

{\bf Proposition 2.4:} {\it En posant
$\D\lambda=\min\Bigl(\lambda_1(Y^{-1});\sqrt{\frac{\pi}{3g}}\Bigr)$, on a la
majoration
$$\int_F\ln f_Y(2;x){\rm d}\nu(x)\le-\frac{\pi}{6\lambda^2}-g\ln\lambda-\frac{g}{2}\ln\frac{6g}{\pi e}\quad.$$}\\

{\it Démonstration:} Soit $t\in]0;2]$. Le lemme 2.3.$\alpha$ implique pour tout
$x\in F$ l'inégalité
$$\ln f_Y(2;x)\le\ln f_Y(t;x)-\pi(2-t)\psi(x)²\quad.$$

Avec les lemmes 2.3.$\beta$, 2.2 et 2.1, on en déduit
$$\begin{array}{rcl}
\D\int_F\ln f_Y(2;x){\rm d}\nu(x)&\le&\D\int_F\ln f_Y(t;x){\rm d}\nu(x)-\pi(2-t)\int_F\psi(x)^2{\rm d}\nu(x)\\
&\le&\D-\frac{g}{2}\ln t-\frac{\pi(2-t)}{12\lambda_1(Y^{-1})^2}\\
\end{array}$$

On obtient le résultat en choisissant $\D t=\frac{6g\lambda^2}{\pi}$.
$\square$\\

\section{Minoration de hauteur}

\subsection{Généralités}

{\bf Définition:} Soit $A$ une variété abélienne complexe de dimension $g\ge1$.
Posons $T_A=\Gamma(A;\Omega_{A/\mathbb{C}})^\vee$ et notons $\Gamma_A$ le réseau
des périodes de $A$ (on a donc un isomorphisme
$A(\mathbb{C})\simeq T_A/\Gamma_A$ de groupes analytiques).

Soit $L:A\rightarrow A^\vee$ une polarisation de $A$. Elle induit une forme de
Riemann $H_L$ sur $T_A$ ({\it i.e.} une forme hermitienne sur $T_A$ telle que
${\rm Im}H_L(\gamma;\delta)\in\mathbb{Z}$ pour tout
$(\gamma;\delta)\in\Gamma_A^2$). Le {\bf diamètre d'injectivité} de $L$ est le
réel $\D\rho(A;L)=\min_{\gamma\in\Gamma_A-\{0\}}\sqrt{H_L(\gamma;\gamma)}$.\\

{\bf Définition:} Soit $(A;L)$ une variété abélienne complexe
principalement polarisée. Désignons par $\nu_1$ la mesure de Haar sur
$A(\mathbb{C})$ de masse 1.

Prenons un faisceau inversible ${\cal L}$ ample sur $A$ définissant $L$ (on a
donc $h^0(A;{\cal L})=1)$, une métrique du cube $\|\ \|$ sur ${\cal L}$
({\it i.e.} une métrique sur $L$ à courbure invariante par translations), et
une section $s\in\Gamma(A;{\cal L})-\{0\}$. On pose
$$I(A;L)=-\int_{A(\mathbb{C})}\ln\|s\|{\rm d}\nu_1+\frac{1}{2}\ln\int_{A(\mathbb{C})}\|s\|^2{\rm d}\nu_1\quad.$$
On vérifie facilement que ce réel ne dépend pas du choix de
$({\cal L};\|\ \|;s)$.\\

Soit $K$ un corps de nombres. Soit $A$ une variété abélienne de dimension $g$
sur $K$.\\

{\bf Définition:} On définit la {\bf hauteur de Faltings} $h_{\rm Fa}(A)$ de $A$
de la manière suivante. Choisissons une extension finie $K'$ de $K$ telle que
$A_{K'}$ soit semi-stable sur $K'$. Désignons par $X$ le modèle de Néron de
$A_{K'}$ sur $B={\rm Spec}(O_{K'})$, par $0_X\in X(B)$ sa section neutre, et par
$\mbox{\Large$\omega$}_X$ le faisceau inversible
$\mbox{\Large$\omega$}_X=0_X^*\Lambda^g\Omega_{X/B}$ sur $B$.

Munissons $\mbox{\Large$\omega$}_X$ de la métrique $\|\ \|_{\rm Fa}$ telle que
pour tout $\sigma\in B(\mathbb{C})=G_{K'}$ et tout
$\varphi\in\Gamma(B_\sigma;\mbox{\Large$\omega$}_{X\sigma})=\Gamma(A_\sigma;\Lambda^g\Omega_{A_\sigma/\mathbb{C}})$, on ait
$\D\|\varphi\|_{\rm Fa}^2(\sigma)=\frac{i^{g^2}}{2^g}\int_{A_\sigma(\mathbb{C})}\varphi\wedge\overline{\varphi}$.

On pose alors $\D h_{\rm Fa}(A)=\frac{\widehat{\deg}(\mbox{\Large$\omega$}_X;\|\ \|_{\rm Fa})}{[K':\mathbb{Q}]}$ (cela ne dépend pas du choix de $K'$).\\

On va utiliser l'inégalité de Bost suivante ({\it cf.} théorème \S 3 de
\cite{Bo2}).\\

{\bf Théorème 3.1 (Bost):} {\it Soit $(A;L)$ une $K$-variété abélienne
de dimension $g$, principalement polarisée. On a alors la minoration
$$h_{\rm Fa}(A)\ge-\frac{g}{2}\ln(2\pi^2)+\frac{2}{[K:\mathbb{Q}]}\sum_{\sigma\in G_K}I(A_\sigma;L_\sigma)\quad.$$}\\

{\it Démonstration:} Voir l'appendice de \cite{GaRe}. $\square$\\

\subsection{Démonstration du théorème 1.1}

Grâce au théorème 3.1, il suffit de minorer $I(A;L)$ en fonction de $\rho(A;L)$
pour toute variété abélienne complexe $(A;L)$ principalement polarisée.

Notons $\mathbb{H}_g$ l'espace de Siegel des matrices
$\Omega\in{\rm M}_g(\mathbb{C})$ symétriques telles que ${\rm Im}\Omega$ soit
définie positive. \`A tout $\Omega\in\mathbb{H}_g$ on associe la fonction thêta
définie par
$$\forall z\in\mathbb{C}^g\quad\theta_\Omega(z)=\sum_{n\in\mathbb{Z}^g}\exp(i\pi ^{\rm t}\!n\Omega n+2i\pi ^{\rm t}\!nz)\quad.$$

Soit $(A;L)$ une $\mathbb{C}$-variété abélienne de dimension $g$,
principalement polarisée. Fixons un $\Omega\in\mathbb{H}_g$ tel que
$A(\mathbb{C})\simeq\mathbb{C}^g/(\mathbb{Z}^g+\Omega\mathbb{Z}^g)$, que $L$
soit induite par $\Theta={\rm div}(\theta_\Omega)$, et que $\Omega$ soit réduite
au sens de Siegel (voir \S V.4 de \cite{Igu}).

En posant $Y={\rm Im}\Omega$, on a en particulier
$\D\lambda_1(Y)^2\ge\frac{\sqrt{3}}{2}$.\\

Désignons par $\nu_1$ la mesure de Haar sur $A(\mathbb{C})$ et par $\nu$ la
mesure de Lebesgue sur $\mathbb{R}^g$. Posons ${\cal L}={\cal O}_A(\Theta)$ et
notons $s$ la section globale de ${\cal L}$ définie par $\Theta$. On munit
${\cal L}$ de la métrique $\|\ \|$ définie par
$$\forall
z=x+iy\in\mathbb{C}^g\quad\|s\|(z)=\sqrt[4]{\det(Y)}\exp(-\pi ^{\rm t}\!yY^{-1}y)|\theta_\Omega(z)|\quad.$$

C'est une métrique du cube sur ${\cal L}$ (voir \S 3 de \cite{Mor}) et on a
$$\ln\int_{A(\mathbb{C})}\|s\|^2{\rm d}\nu_1=-\frac{g}{2}\ln2\quad.$$

Maintenant, majorons le terme $\D\int_{A(\mathbb{C})}\ln\|s\|{\rm d}\nu_1$.
Posons $F=[0;1]^g$. En utilisant l'inégalité de Jensen puis la formule de
Parseval, on obtient pour tout $y\in F$:
$$\begin{array}{rcl}
\D\int_F\ln\|s\|^2(x+\Omega y){\rm d}\nu(x)&\le&\D\ln\int_F\|s\|^2(x+\Omega y){\rm d}\nu(x)\\
&=&\D\ln\Bigl[\sqrt{\det(Y)}\sum_{n\in\mathbb{Z}^g}\exp[-2\pi ^{\rm t}\!(y+n)Y(y+n)]\Bigr]\\
&=&\D\ln f_Y(2;y)\quad.\\
\end{array}$$

On pose
$\D\lambda=\min\Bigl(\lambda_1(Y^{-1});\sqrt{\frac{\pi}{3g}}\Bigr)$. \`A l'aide
de la proposition 2.4, on en déduit
$$\begin{array}{rcl}
\D2I(A;L)&\ge&\D-\frac{g}{2}\ln2-\int_F\ln f_Y(2;y){\rm d}\nu(y)\\
&\ge&\D\frac{\pi}{6\lambda^2}+g\ln\lambda+\frac{g}{2}\ln\frac{3g}{\pi e}\quad.\\
\end{array}$$

On conclut en appliquant le lemme 3.2 ci-dessous. $\square$\\

{\bf Lemme 3.2:} {\it En posant
$\D\rho=\min\Bigl(\rho(A;L);\sqrt{\frac{\pi}{3g}}\Bigr)$, on a l'égalité
$\rho=\lambda$.}\\

{\it Démonstration:} Identifions $T_A$ à $\mathbb{C}^g$ et $\Gamma_A$ à
$\mathbb{Z}^g+\Omega\mathbb{Z}^g$. La forme de Riemann $H_L$ vérifie alors
$$\forall \gamma=m+\Omega n\in\Gamma_A\quad H_L(\gamma;\gamma)={}^{\rm t}\!\bar{\gamma}Y^{-1}\gamma={}^{\rm t}\!(m+Xn)Y^{-1}(m+Xn)+{}^{\rm t}\!nYn\quad.$$

Si $g=1$, on vérifie aisément que
$\D\rho(A;L)=\frac{1}{\sqrt{Y}}=\lambda_1(Y^{-1})$.\\

Supposons maintenant $g\ge2$. Soit $\gamma=m+\Omega n\in\Gamma_A-\{0\}$.
Si $n\neq0$ on a $\D H(\gamma;\gamma)\ge{}^{\rm t}\!nYn\ge\frac{\sqrt{3}}{2}\ge\frac{\pi}{3g}$, et si $n=0$ on a
$H(\gamma;\gamma)={}^{\rm t}\!mY^{-1}m\ge\lambda_1(Y^{-1})^2$.
D'où $\rho\ge\lambda$.

De même, pour tout $m'\in\mathbb{Z}^g-\{0\}$, on a
$|\!|m'|\!|_{Y^{-1}}^2=H(m';m')\ge\rho(A;L)^2$. Donc $\lambda\ge\rho$.
$\square$\\

\ \\

{\small Pascal Autissier. I.M.B., université Bordeaux I, 351, cours de la
Libération, 33405 Talence cedex, France.

pascal.autissier@math.u-bordeaux1.fr}

\end{document}